\documentclass[11pt,draft]{amsart}
\usepackage[shortalphabetic]{amsrefs}
\usepackage{amsfonts}
\usepackage{amsmath}
\usepackage{amsthm}
\usepackage{pstricks}
\usepackage{pst-node}
\usepackage{pst-plot}

\def\Ee{\mathbb{E}}
\def\nr{2pt}
\def\O{\mathcal{O}}
\def\Pp{\mathbb{P}}
\def\Qq{\mathbb{Q}}

\newcommand{\M}[2]{\mathcal{M}_{#1,#2}}
\newcommand{\Mct}[2]{\mathcal{M}^{\text{ct}}_{#1,#2}}
\newcommand{\Mrt}[2]{\mathcal{M}^{\text{rt}}_{#1,#2}}
\newcommand{\Mla}[3]{\mathcal{M}^{\lambda_#1}_{#2,#3}}
\newcommand{\Mch}[3]{\mathcal{M}^{\ch_{#1}}_{#2,#3}}
\newcommand{\Mbar}[2]{\overline{\mathcal{M}}_{#1,#2}}
\newcommand{\im}{\operatorname{im}}
\newcommand{\ch}{\operatorname{ch}}
\newcommand{\irr}{\operatorname{irr}}

\newcommand{\iibar}{\overline{2i}}
\newcommand{\nbar}{\overline{n}}

\newcommand{\smallsetminus}{\mbox{\footnotesize{$\setminus$}}}

\def\gradius{.15cm}
\newcommand{\genusvertex}[3]{\rput#1{\cnode{\gradius}{#2}}\rput#1{\scriptsize $#3$}}
\newcommand{\markedpoint}[3]{\rput#1{\rnode{#2}{\scriptsize $#3$}}}

\theoremstyle{plain}
\newtheorem{lemma}[subsection]{Lemma}
\newtheorem{proposition}{Proposition}
\newtheorem{corollary}[proposition]{Corollary}
\newtheorem*{startheorem}{Theorem $\star$}
\newtheorem*{question}{Question}

\theoremstyle{definition}
\newtheorem{definition}[subsection]{Definition}
\newtheorem{remark}{Remark}

\title{Tautological pairings on moduli spaces of curves}

\author{Renzo Cavalieri}
\address{Colorado State University, Department of Mathematics, Weber
  Building, Fort Collins, CO 80523-1874}
\email{renzo@math.colostate.edu}

\author{Stephanie Yang}
\address{Institutionen f\"or Matematik, Kungliga Tekniska
  H\"ogskolan, 100 44 Stockholm, Sweden}
\email{stpyang@math.kth.se}

\begin{document}
\psset{unit=.75cm}

\begin{abstract}
We discuss analogs of Faber's conjecture for two nested sequences of
partial compactifications of the moduli space of smooth curves.  We
show that their tautological rings are one-dimensional in top degree
but do not satisfy Poincar\'e duality.
\end{abstract}

\maketitle

The structure of the tautological ring of the moduli space of stable
curves is predicted by the {\em Faber conjecture}, which states that
$R^\ast(\Mbar{g}{n})$ is Gorenstein with socle in dimension
$3g-3+n=\dim\Mbar{g}{n}$. We break this statement into two parts:

\begin{description}
\item[Socle] The tautological ring vanishes in high degree and
it is one-dimensional in top degree.
\begin{equation} R^k(\Mbar{g}{n}) = 
\begin{cases}
 0,&\text{ if $k>3g-3+n$} \\
 \Qq,&\text{ if $k=3g-3+n$}
\end{cases}
\end{equation}
\item[Poincar\'e duality] For $0\leq k\leq 3g-3+n$, the bilinear
  pairing
\begin{align}
R^k(\Mbar{g}{n})\times R^{3g-3+n-k}(\Mbar{g}{n})&\to R^{3g-3+n}(\Mbar{g}{n})
\end{align}
is non-degenerate.
\end{description}

In \cite{MR1786488}, Faber and Pandharipande speculate that the
tautological rings of $\Mct{g}{n}$ and $\Mrt{g}{n}$, two partial
compactifications of the moduli space of curves, satisfy analogous
properties. While the socle statements have since been proven in all
instances, Poincar\'e duality remains open. We give evidence
that the two properties are not necessarily immediately correlated.
  
We define two chains of partial compactifications (Defs.~\ref{def:mla}
and~\ref{def:mch})
\begin{gather}
  \Mct{g}{n}=\Mla{g}{g}{n}\subseteq\Mla{g-1}{g}{n}\subseteq\cdots\subseteq\Mla{1}{g}{n}\subseteq\Mbar{g}{n}\\
  \Mrt{g}{n}=\Mch{2g-1}{g}{n}\subseteq\Mch{2g-3}{g}{n}\subseteq\cdots\subseteq\Mch{1}{g}{n}\subseteq\Mbar{g}{n}
\end{gather}
and define their tautological rings by restriction.  The main results
of this paper address the analog of Faber's conjecture for these
spaces.  The socle statement extends but Poincar\'e duality fails.

\begin{proposition}\label{prop:socle}
\label{socle}
For $i=1,2,\ldots,g$,
\begin{enumerate}
\item\label{prop:mla} $R^k(\Mla{i}{g}{n})=0$ for $k>3g-3+n-i$ and
    $R^{3g-3+n-i}(\Mla{i}{g}{n})\cong\Qq$,
\item\label{prop:mch} $R^k(\Mch{2i-1}{g}{n})=0$ for $k>3g-2+n-2i$ and
    $R^{3g-2+n-2i}(\Mch{2i-1}{g}{n})\cong\Qq$.
\end{enumerate}
\end{proposition}

\begin{proposition}\label{prop:duality}
  For $g\geq 2$ and $(g,n)\neq(2,0)$, the pairing
\begin{equation}
  R^g(\Mla{1}{g}{n})\times R^{2g-4+n}(\Mla{1}{g}{n})\to R^{3g-4+n}(\Mla{1}{g}{n})
\end{equation}
is not perfect.
\end{proposition}

\begin{proposition}\label{prop:lai}
For any fixed $g$ and $2\leq i \leq g-1$, and for every $n\gg 0$,
either the tautological restriction sequence
\begin{equation}\label{eq:restrict}
  R^{i-1}(\Mbar{g}{n}\smallsetminus\Mla{i}{g}{n})\longrightarrow R^g(\Mbar{g}{n})\longrightarrow R^g(\Mla{i}{g}{n})\longrightarrow 0
\end{equation}
is not exact in the middle, or the pairing
\begin{equation}
\label{eq:pairi}
  R^g(\Mla{i}{g}{n})\times R^{2g-3-i+n}(\Mla{i}{g}{n})\to R^{3g-3-i+n}(\Mla{i}{g}{n})
\end{equation}
is not perfect.  
\end{proposition}

The Chow ring is known to be tautological in codimensions $0$ and $1$
(\cites{MR717614}). Therefore Proposition
\ref{prop:lai} immediately implies:

\begin{corollary}\label{cor:lai}
  If $g\geq 3$ and $i=2$, then the pairing
  (\ref{eq:pairi}) is not perfect for every $n\gg 0$.
\end{corollary}

We also show that for $g\geq 2$, the dimension of the kernel of the
map
\begin{equation}\label{eq:phifirst}
  R^g(\Mla{1}{g}{n})\to\hom\left(R^{2g-4+n}(\Mla{1}{g}{n}),R^{3g-4+n}(\Mla{1}{g}{n})\right).
\end{equation}
becomes arbitrarily large as either $g$ or $n$
grows (Corollary~\ref{cor:dimension}).

In Lemma~\ref{le:extremal} we show that for the extremal case $i=g$,
$\Mla{g}{g}{n}=\Mct{g}{n}$ and $\Mch{2g-1}{g}{n}=\Mrt{g}{n}$. Our
results therefore answer extensions of speculations about these spaces
that appeared in \cite{MR1786488}.  The socle statements of the Faber
conjectures for $\Mrt{g}{n}$, $\Mct{g}{n}$, and $\Mbar{g}{n}$ were
proved by Graber-Vakil (\cite{MR2176546}*{\S 5.5--\S 5.7}) and
Faber-Pandharipande (\cite{MR2120989}*{\S4.1}).  Prior to these
results, the socle statement for $\mathcal{M}_g$ was shown by
Looijenga and Faber
(\citelist{\cite{MR1346214}\cite{MR1722541}*{Theorem 2}}).  The
Poincar\'e duality property is only known for $g=0$ by Keel
(\cite{MR1034665}), and for $\mathcal{M}_g$ for $g\leq 23$ by
Faber~(\cite{MR1722541}).

\begin{remark}
The tautological restriction sequence for $M_g$ is exact in degrees
$\geq g-1$, and exactness is conjectured in all degrees for the
$\Mrt{g}{n}$ and $\Mct{g}{n}$ (\cite{MR2120989}).
\end{remark}

\begin{remark}
It is easy to see that Poincare' duality fails in arbitrary degrees by
taking nonzero elements of high degree in the ideal generated by the
kernel of the map~\eqref{eq:phifirst}.
\end{remark}

\begin{remark}
Note that $R^\ast(\Mla{1}{2}{0})$ is Gorenstein. The socle dimension
is 2 and it is straightforward to check that the intersection matrix
for the two generators of $R^1(\Mla{1}{2}{0})$ is non-degenerate.
Since $\pi_1^\ast(\Mla{1}{2}{0})=\Mla{1}{2}{1}$,
Proposition~\ref{prop:duality} shows that if $\mathcal{M}$ is a moduli
space which satisfies Faber's conjecture, its universal family does
not necessarily satisfy Faber's conjecture.
\end{remark}


%
\begin{remark}
The classes $\lambda_i$ and $\ch_{2i-1}$ vanish respectively on
$\Mbar{g}{n}\smallsetminus\Mla{i}{g}{n}$ and
$\Mbar{g}{n}\smallsetminus\Mch{2i-1}{g}{n}$.  This motivates our
notation.
\end{remark}

This paper is organized as follows: in \S\ref{sec:background} we
recall some basic definitions and provide references to the existing
literature. In each of the subsequent sections we prove the three
propositions.

The authors thank Carel Faber and Ravi Vakil for helpful
discussions.  The research was partially supported the National
Science Foundation, the Gustafsson Foundation, and the Wallenberg
Foundation.

\section{Background}\label{sec:background}

 The {\em tautological ring} $R^\ast(\Mbar{g}{n})$ is a
natural subring of the Chow ring $A^\ast(\Mbar{g}{n})$ 
elegantly defined 
in~\cite{MR2120989}: as $g$ and $n$ vary, the tautological rings form
the smallest system of $\Qq$-subalgebras of $A^\ast(\Mbar{g}{n})$ that
are closed under the natural forgetful morphisms
\begin{equation}
  \pi_{i}\colon\Mbar{g}{n+1}\to\Mbar{g}{n},\\
\end{equation}
and the gluing morphisms
\begin{gather}
  \iota_{\irr}\colon\Mbar{g-1}{n+2}\to\Mbar{g}{n} \\
  \iota_{g_1,n_1}\colon\Mbar{g_1}{n_1+1}\times\Mbar{g_2}{n_2+1}\to\Mbar{g_1+g_2}{n_1+n_2}.
\end{gather}
The tautological ring contains boundary strata $\delta_\Gamma$ (the
closure of the locus of curves whose dual graph is $\Gamma$),
cotangent line classes $\psi_i$, Mumford-Morita $\kappa$ classes,
chern classes $\lambda_i$ of the Hodge bundle.

Much is known about the intersection theory of such classes. An
excellent, albeit unfinished and unpublished, reference is
\cite{jknotes}. Other references include \cite{MR1644323},
\cite{MR2003030}, and \cite{MR717614}.

The following formula will be used in the proof of
Proposition~\ref{prop:duality}.
\begin{lemma}\label{le:kappakappa}
For any value of $n$ for which the integrals are defined:
\begin{align}
	\int_{\Mbar{0}{n}}\kappa_{n-3} &=
        1 \label{a}\\ \int_{\Mbar{1}{n}}\kappa_{n} &=
        \frac{1}{24} \label{b}\\ \int_{\Mbar{0}{n}}\kappa_i\kappa_{n-i}
        &= \binom{n-1}{i+1}-1 \label{c}
  \end{align}
\end{lemma}
\begin{proof}
Equations (\ref{a}) and (\ref{b}) follow immediately from the pullback
formula for $\kappa$ classes,~(\cite{MR1486986}*{Eq.~(1.10)}),
\begin{equation}
  \kappa_a=\pi_{n+1}^\ast(\kappa_a)-\psi_{n+1}^a.
\end{equation}
\begin{align}
  \int_{\Mbar{0}{n}} \kappa_i\kappa_{n-i-3}
  &=\int_{\Mbar{0}{n}}\left(\pi_*\psi_{n+1}^{i+1}\right)\kappa_{n-i-3} \\
  &=\int_{\Mbar{0}{n+1}}\psi_{n+1}^{i+1}\left(\kappa_{n-i-3}-\psi_{n+1}^{n-i-3}\right)\label{eq:ac}\\
  &=\int_{\Mbar{0}{n+1}}\psi_{n+1}^{i+1}\left(\pi_*\psi_{n+2}^{n-i-2}\right)-\int_{\Mbar{0}{n+1}}\psi_{n+1}^{n-2}\\
  &=\int_{\Mbar{0}{n+2}}\left(\psi_{n+1}^{i+1}-D_{n+1,n+2}\right)\psi_{n+2}^{n-i-2}-\int_{\Mbar{0}{n+1}}\psi_{n+1}^{n-2}\label{eq:rv}\\
  &=\int_{\Mbar{0}{n+2}}\psi_{n+1}^{i+1}\psi_{n+2}^{n-i-2}-0-\int_{\Mbar{0}{n+1}}\psi_{n+1}^{n-2}\\
  &=\binom{n-1}{i+1}-1.
\end{align}
\end{proof}

\begin{definition}\label{def:mla}
$\Mla{i}{g}{n}\subseteq\Mbar{g}{n}$ is the locus of curves whose dual
graph has genus $\leq g-i$. Equivalently, $\Mla{i}{g}{n}$ is the locus
of curves where the sum of the geometric genera of the components is
at least $i$.
\end{definition}
\begin{definition}\label{def:mch}
$\Mch{2i-1}{g}{n}\subseteq\Mbar{g}{n}$ is the locus of curves with
at least one component of genus at least $i$.
\end{definition}

\begin{lemma}\label{le:extremal}
\label{ext}
  $\Mla{g}{g}{n}=\Mct{g}{n}$ and $\Mch{2g-1}{g}{n}=\Mrt{g}{n}$.
\end{lemma}
\begin{proof}
  The dual graph of any curve in $\Mla{g}{g}{n}$ is connected of genus
  $0$ and thus is a tree.  Any curve in $\Mch{2g-1}{g}{n}$ has at
  least one (hence exactly one) component of genus $g$; the other
  components must necessarily form  trees of rational curves.
\end{proof}

\begin{lemma}\label{le:vanish}
  For $i=1,2,\ldots,g$,
  \begin{enumerate}
  \item The class $\lambda_i$ vanishes on
    $\Mbar{g}{n}\smallsetminus\Mla{i}{g}{n}$.
  \item The class $\ch_{2i-1}$ vanishes on
    $\Mbar{g}{n}\smallsetminus\Mch{2i-1}{g}{n}$
  \end{enumerate}
\end{lemma}
\begin{proof}

The boundary stratum $\delta_\Gamma$ is the image of the gluing
morphism
\begin{equation}
  \delta_\Gamma =
  \im\left(\iota_\Gamma\colon\prod{\Mbar{g_j}{n_j+d_j}}\to\Mbar{g}{n}\right).
\end{equation}
The Hodge bundle splits when restricted to $\delta_\Gamma$
as in~\cite{MR1728879}*{Eqs.~(17)~and~(18)}:
\begin{equation}
  \iota_{\Gamma}^\ast(\Ee) = \oplus \Ee_{g_j,n_j+d_j} \oplus \O^n.
\end{equation}

To see (a), use the Whitney formula,
\begin{equation}
  \label{w}
 \iota_\Gamma^\ast( c(\Ee)) = c(\iota_\Gamma^\ast\Ee) =
  \prod\left(1+\lambda_{1,j}+\cdots+\lambda_{g_j,j}\right)
\end{equation}
where $\lambda_{i,j}$ is the $i$-th chern class of the Hodge bundle of
the $j$-th factor $\Mbar{g_j}{n_j+d_j}$. If
$\delta_\Gamma\in\Mbar{g}{n}\smallsetminus\Mla{i}{g}{n}$, then $\sum
g_j<i$ and term of degree $i$ in equation~\eqref{w} vanishes.

To see (b), use the additivity of the chern character, on
$\delta_\Gamma$:
\begin{equation}
 \iota_\Gamma^\ast(\ch_{2i-1}(\Ee))= \ch_{2i-1}(\iota_\Gamma^\ast\Ee) =
  \sum_j \ch_{2i-1,j}
\end{equation}
where $\ch_{2i-1,j}$ is the $(2i-1)$-th chern character of the Hodge
bundle of the $j$-th factor $\Mbar{g_j}{n_j+d_j}$. Since
$\delta_\Gamma\in\Mbar{g}{n}\smallsetminus\Mch{2i-1}{g}{n}$, all $g_j<i$
and thus $\ch_{2i-1}=0$.
\end{proof}


We conclude this section by recalling theorem $\star$ by Graber-Vakil
(\cite{MR2176546}), which is a key ingredient in the proof of
Proposition~\ref{prop:socle}.

\begin{startheorem} Any tautological class of degree $k$ on $\M{g}{n}$
  vanishes on the open set consisting of strata satisfying
\begin{equation}
\#\text{\rm\ genus 0 components} < k-g+1.
\end{equation}
\end{startheorem}

\section{Socle}

\label{soclesec}
In this section we prove Proposition~\ref{socle}. The strategy of the
proof is natural: theorem $\star$ forces tautological classes of high
degree to be supported on strata with many rational components. On the
other hand, curves in $\Mla{i}{g}{n}$ and $\Mch{2i-1}{g}{n}$ satisfy
geometric conditions that limit the number of rational
components. These constraints imply high degree vanishing and force
tautological classes in the socle degree to be supported on exactly
one boundary stratum up to rational equivalence.

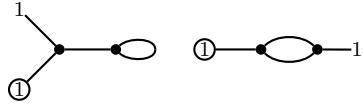
\begin{figure}[tb]\label{fig:feynman}
\centerline{
\begin{pspicture}(-.707,-.707)(1.7,0.707)
\rput(0,0){
  \markedpoint{(-.707,0.707)}{mp1}{1}
  \genusvertex{(-.707,-.707)}{v1}{1}
  \cnode*(0,0){2pt}{v2}
  \cnode*(1,0){2pt}{v3}
  \pnode(1.7,0){x4}
  \ncline{mp1}{v2}
  \ncline{v1}{v2}
  \ncline{v2}{v3}
  \nccurve[angleA=45,angleB=90]{v3}{x4}
  \nccurve[angleA=-45,angleB=-90]{v3}{x4}
}
\end{pspicture}
\begin{pspicture}(-.707,-.707)(1.7,0.707)
\rput(0,0){
  \genusvertex{(0,0)}{v1}{1}
  \cnode*(1,0){2pt}{v2}
  \cnode*(2,0){2pt}{v3}
  \markedpoint{(2.7,0)}{mp1}{1}
  \pnode(1.7,0){x4}
  \ncline{v1}{v2}
  \nccurve[angleA=45,angleB=135]{v2}{v3}
  \nccurve[angleA=-45,angleB=-135]{v2}{v3}
  \ncline{v3}{mp1}
}
\end{pspicture}
}
\caption{These two dual graphs differ by a Feynman move.}
\end{figure}

\begin{definition}
  A Feynman move replaces a portion of a graph of type (a) on a dual graph with
  one of type (b) or (c), as illustrated below:
  
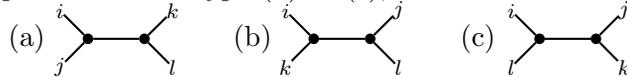
\begin{center}
  \begin{pspicture}(-.5,-.5)(10.5,0.5)
    \rput(0,0){(a)}
    \rput(1,0){
      \cnode*(0,0){2pt}{v1}
      \cnode*(1,0){2pt}{v2}
      \markedpoint{(-.5,0.5)}{m1}{i}
      \markedpoint{(-.5,-.5)}{m2}{j}
      \markedpoint{(1.5,0.5)}{m3}{k}
      \markedpoint{(1.5,-.5)}{m4}{l}
      \ncline{v1}{v2}
      \ncline{v1}{m1}
      \ncline{v1}{m2}
      \ncline{v2}{m3}
      \ncline{v2}{m4}
    }
    \rput(4,0){(b)}
    \rput(5,0){
      \cnode*(0,0){2pt}{v1}
      \cnode*(1,0){2pt}{v2}
      \markedpoint{(-.5,0.5)}{m1}{i}
      \markedpoint{(-.5,-.5)}{m3}{k}
      \markedpoint{(1.5,0.5)}{m2}{j}
      \markedpoint{(1.5,-.5)}{m4}{l}
      \ncline{v1}{v2}
      \ncline{v1}{m1}
      \ncline{v1}{m3}
      \ncline{v2}{m2}
      \ncline{v2}{m4}
    }
    \rput(8,0){(c)}
    \rput(9,0){
      \cnode*(0,0){2pt}{v1}
      \cnode*(1,0){2pt}{v2}
      \markedpoint{(-.5,0.5)}{m1}{i}
      \markedpoint{(-.5,-.5)}{m4}{l}
      \markedpoint{(1.5,0.5)}{m2}{j}
      \markedpoint{(1.5,-.5)}{m3}{k}
      \ncline{v1}{v2}
      \ncline{v1}{m1}
      \ncline{v1}{m4}
      \ncline{v2}{m2}
      \ncline{v2}{m3}
    }
    \end{pspicture}
  \end{center}
  The half edges $i$, $j$, $k$, and $l$ may be
  glued to other half edges to form edges (See Figure \ref{fig:feynman}).
\end{definition} 

\begin{lemma}
  If the dual graphs of two boundary strata differ by Feynman moves,
  then they are rationally equivalent.
\end{lemma}
\begin{proof}
  This is immediate by noting that $\Mbar{0}{4}\cong\Pp_1$, so its
  boundary points are rationally equivalent.  This equivalence is
  preserved under gluing morphisms.
\end{proof}

\begin{remark}\label{rem:feynman}
It is a standard combinatorial fact that two trivalent graphs with
same number of vertices and edges differ by a finite number of
Feynman moves (\cite{MR0002545}).
\end{remark}

The proof of Proposition \ref{socle} now follows from some careful
bookkeeping.

\begin{proof}[Part~\ref{prop:mla}]
  Theorem $\star$ implies that any class in
  $R^{3g-3+n-i+k}(\Mla{i}{g}{n})$ $(k\geq 0)$ must be supported on
  boundary strata containing at least $2g+n-2-i+k$ rational
  components. Let $\delta_\Gamma$ be any one such boundary stratum. By
  the definition of $\Mla{i}{g}{n}$, the sum of the genera of the
  vertices of graph $\Gamma$ is at most $g-i$.

  Stability implies that the incidence of any rational vertex is at
  least three, and there are at most $n$ marked half-edges incident to
  these vertices.  Thus there are at least
  \begin{equation}  
    3(2g+n-2-i+k)-n=6g+2n-6+3i+3k
  \end{equation}
  half edges which must be glued to other half-edges.  Since $g_\Gamma
  \leq g-i$, at least
  \begin{equation}
    6g+2n-6-3i+3k-2(3g+n-2-2i+k)=i+k
  \end{equation}
  half-edges must be glued to vertices of positive genus.  The total
  genus of the curve represented by $\Gamma$ is $g$, so the only
  possibility is $k=0$, i.e.~all rational vertices trivalent, and the
  remaining $i$ half edges glued to $i$ vertices of genus $1$ and
  incidence $1$.  By Remark~\ref{rem:feynman}, any two such graphs
  differ by a finite number of Feynman moves, and hence represent
  rationally equivalent tautological classes. It is immediate to check
  that such graphs actually live in degree $3g-3+n-i$.
\end{proof}

\begin{proof}[Part~\ref{prop:mch}] 
  This is similar to the previous proof, except now we require the
  graph $\Gamma$ to have at least one vertex of genus $\geq i$.
  Theorem $\star$ forces any class in $R^{3g+n-2-2i+k}(\Mch{i}{g}{n})$
  ($k\geq 0$) to be supported on boundary strata with at least
  $2g+n-1-2i+k$ rational components, and stability implies that there
  is at least three times as many half edges, only $n$ of which are
  not glued to some other half edge. The total genus of the graph must
  be $g$, the only consistent possibility is for $k=0$, i.e.~all
  rational vertices are trivalent and there is exactly one vertex of
  genus $i$ and incidence $1$. Again, all such strata represent
  rationally equivalent classes because any two graphs are equivalent 
  up to a finite number of Feynman moves.

  The degree of such strata is $3g+n-1-3i$, which means that one must
  have a class of degree $i-1$ decorating the genus $i$ vertex. Since
  $R^{i-1}(\Mrt{i}{1})\cong\Qq$ (this is the socle statement of
  Faber's conjecture for $\Mrt{i}{1}$), there is only one such nonzero
  class up to scalar multiple.
\end{proof}

\begin{figure}[tb]
\centerline{
\begin{pspicture}(-1,-1)(13,1)
\genusvertex{(0,1)}{v0}{1}
\cnode*(1,1){2pt}{v1}
\cnode*(2,1){2pt}{v2}
\rput(3,1){\rnode{v3}{$\cdots$}}
\cnode*(4,1){2pt}{v4}
\cnode*(5,1){2pt}{v5}
\cnode*(6,1){2pt}{v6}
\rput(7,1){\rnode{v7}{$\cdots$}}
\cnode*(8,1){2pt}{v8}
\cnode*(9,1){2pt}{v9}
\cnode*(10,1){2pt}{v10}
\rput(11,1){\rnode{v11}{$\cdots$}}
\cnode*(12,1){2pt}{v12}
\markedpoint{(13,1)}{v13}{n}
\genusvertex{(1,0)}{w1}{1}
\genusvertex{(2,0)}{w2}{1}
\genusvertex{(4,0)}{w4}{1}
\cnode*(5,0){2pt}{w5}
\cnode*(6,0){2pt}{w6}
\cnode*(8,0){2pt}{w8}
\markedpoint{(9,0)}{w9}{1}
\markedpoint{(10,0)}{w10}{2}
\markedpoint{(12,0)}{w12}{n-1}
\pnode(5,-.7){x5}
\pnode(6,-.7){x6}
\pnode(8,-.7){x8}
\ncline{v0}{v1}
\ncline{v1}{v2}
\ncline{v2}{v3}
\ncline{v3}{v4}
\ncline{v4}{v5}
\ncline{v5}{v6}
\ncline{v6}{v7}
\ncline{v7}{v8}
\ncline{v8}{v9}
\ncline{v9}{v10}
\ncline{v10}{v11}
\ncline{v11}{v12}
\ncline{v12}{v13}
\ncline{v1}{w1}
\ncline{v2}{w2}
\ncline{v4}{w4}
\ncline{v5}{w5}
\ncline{v6}{w6}
\ncline{v8}{w8}
\ncline{v9}{w9}
\ncline{v10}{w10}
\ncline{v12}{w12}
\nccurve[angleA=-135,angleB=180]{w5}{x5}\nccurve[angleA=-45,angleB=0]{w5}{x5}
\nccurve[angleA=-135,angleB=180]{w6}{x6}\nccurve[angleA=-45,angleB=0]{w6}{x6}
\nccurve[angleA=-135,angleB=180]{w8}{x8}\nccurve[angleA=-45,angleB=0]{w8}{x8}
\end{pspicture}
}
\bigskip
\centerline{
\begin{pspicture}(-1,-1)(9,1)
\genusvertex{(0,1)}{v0}{i}
\rput(0,1.5){\small $\psi^{i-1}$}
\cnode*(1,1){2pt}{v1}
\cnode*(2,1){2pt}{v2}
\rput(3,1){\rnode{v3}{$\cdots$}}
\cnode*(4,1){2pt}{v4}
\cnode*(5,1){2pt}{v5}
\cnode*(6,1){2pt}{v6}
\rput(7,1){\rnode{v7}{$\cdots$}}
\cnode*(8,1){2pt}{v8}
\markedpoint{(9,1)}{v9}{n}
\cnode*(1,0){2pt}{w1}
\cnode*(2,0){2pt}{w2}
\cnode*(4,0){2pt}{w4}
\markedpoint{(5,0)}{w5}{1}
\markedpoint{(6,0)}{w6}{2}
\markedpoint{(8,0)}{w8}{n-1}
\pnode(1,-.7){x1}
\pnode(2,-.7){x2}
\pnode(4,-.7){x4}
\ncline{v0}{v1}
\ncline{v1}{v2}
\ncline{v2}{v3}
\ncline{v3}{v4}
\ncline{v4}{v5}
\ncline{v5}{v6}
\ncline{v6}{v7}
\ncline{v7}{v8}
\ncline{v8}{v9}
\ncline{v1}{w1}
\ncline{v2}{w2}
\ncline{v4}{w4}
\ncline{v5}{w5}
\ncline{v6}{w6}
\ncline{v8}{w8}
\nccurve[angleA=-135,angleB=180]{w1}{x1}\nccurve[angleA=-45,angleB=0]{w1}{x1}
\nccurve[angleA=-135,angleB=180]{w2}{x2}\nccurve[angleA=-45,angleB=0]{w2}{x2}
\nccurve[angleA=-135,angleB=180]{w4}{x4}\nccurve[angleA=-45,angleB=0]{w4}{x4}
\end{pspicture}
}
\caption{Generators in socle degree for $\Mla{i}{g}{n}$ and
  $\Mch{i}{g}{n}$.}
\end{figure}

\begin{remark}\label{rem:eval}
  It is important to note that $\lambda_i$ (resp. $\ch_{2i-1}$) does
  not vanish on the unique generator of $R^{top}(\Mla{i}{g}{n})$
  (resp. $R^{top}(\Mch{2i-1}{g}{n})$), and hence it can be used as an
  evaluation class: multiplication by $\lambda_i$ gives an isomorphism
  between the socle and $R^{3g-3+n}(\overline{\mathcal{M}}_{g,n})$.
\end{remark}

\section{Failure of Poincar\'e duality for $\Mla{1}{g}{n}$}
\label{failure}

\begin{figure}[tb]
\begin{center}
\begin{pspicture}(-1,-1)(6,1)
\rput(.9,-.6){\small $\Gamma(0,1,1,\emptyset)$}
\rput(0,0){
  \genusvertex{(0,0)}{v1}{1}
  \cnode*(1,0){\nr}{v2}
  \pnode(1.7,0){x2}
  \markedpoint{(1,0.5)}{p1}{}
  \ncline{v1}{v2}
  \ncline{v2}{p1}
  \nccurve[angleA=45,angleB=90]{v2}{x2}
  \nccurve[angleA=-45,angleB=-90]{v2}{x2}  
}
\rput(3.9,-.6){\small $\Gamma(0,1,1,\{1\})$}
\rput(3,0){
  \genusvertex{(0,0)}{v2}{1}
  \cnode*(1,0){\nr}{v1}
  \pnode(1.7,0){x1}
  \markedpoint{(0,0.5)}{p1}{}
  \ncline{v1}{v2}
  \ncline{v2}{p1}
  \nccurve[angleA=45,angleB=90]{v1}{x1}
  \nccurve[angleA=-45,angleB=-90]{v1}{x1}  
}
\end{pspicture}
\\
\begin{pspicture}(-1,-1)(4,1)
\rput(0.5,-.5){\small $A$}
\rput(0,0){
  \cnode*(0,0){\nr}{v1}
  \cnode*(1,0){\nr}{v2}
  \pnode(-.7,0){x1}
  \pnode(1.7,0){x2}
  \markedpoint{(1,0.5)}{p1}{}
  \ncline{v1}{v2}
  \ncline{v2}{p1}
  \nccurve[angleA=135,angleB=90]{v1}{x1}
  \nccurve[angleA=-135,angleB=-90]{v1}{x1}  
  \nccurve[angleA=45,angleB=90]{v2}{x2}
  \nccurve[angleA=-45,angleB=-90]{v2}{x2}  
}
\rput(3,-.5){\small $B$}
\rput(3,0){
  \cnode*(0,0){\nr}{v1}
  \pnode(-.7,0){x1}
  \pnode(0.7,0){x2}
  \markedpoint{(0,0.5)}{p1}{}
  \ncline{v1}{p1}
  \nccurve[angleA=135,angleB=90]{v1}{x1}
  \nccurve[angleA=-135,angleB=-90]{v1}{x1}  
  \nccurve[angleA=45,angleB=90]{v1}{x2}
  \nccurve[angleA=-45,angleB=-90]{v1}{x2}  
}
\end{pspicture}
\end{center}
\caption{The graphs $\Gamma(0,1,1,\emptyset)$ and
$\Gamma(0,1,1,\{1\})$ for $(g,n)=(2,1)$, and the corresponding graphs
$A$ and $B$.}\label{fig:Rpair}
\end{figure}
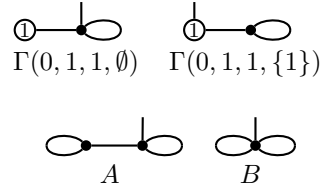

In this section we construct counterexamples to the Poincar\'e duality
part of Faber's conjectures for $\Mla{1}{g}{n}=\Mch{1}{g}{n}$ for
$g\geq 2$ and $(g,n)\neq(2,0)$.  Choose any triple $(a,b,c)$ of
integers satisfying $a+b+c=g$, where $a$ is non-negative and $b$ and
$c$ are positive.  Choose any subset $S$ of the $n$ points.  Let
$\Gamma(a,b,c,S) \in R^g(\Mbar{g}{n})$ denote the graph with two
vertices connected by $b$ edges: one genus 1 vertex with $a$
self-edges and carrying the points in $S$, and one genus 0 vertex with
$c$ self-edges and carrying the points in $S^c$ (see
Figure~\ref{fig:Rpair}). Let $\delta(a,b,c,S)$ denote the associated
boundary stratum.

When $(a,c,S)\neq(c-1,a+1,S^c)$, the strata $\delta(a,b,c,S)$ and
$\delta(c-1,b,a+1,S^c)$ are not rationally equivalent, but their
difference lies in the kernel of multiplication by $\lambda_1$.

\begin{proof}[Proof of Proposition~\ref{prop:duality}]
Let $\gamma_1=\delta(a,b,c,S)$ and
$\gamma_2=\delta(c-1,b,a+1,S^c)$. Note that if $(a,c,S)=(c-1,a+1,S^c)$
then $\gamma_1$ and $\gamma_2$ are the same. Therefore assume that
this is not the case; we can always do this if $g\geq 2$ and
$(g,n)\neq(2,0)$. The fact that $\gamma_1-\gamma_2$ lies in the kernel
of the map $\phi$ from~\eqref{eq:phi} follows from the fact that
$\lambda_1$ vanishes on $\Mbar{0}{k}$ and $\lambda_1$ is equivalent to
$\frac{1}{24}\delta_{\irr}$ on $\Mbar{1}{k}$, for all $k$ for which
these spaces are defined.  Thus $\lambda_1\cdot\gamma_1$ and
$\lambda_1\cdot \gamma_2$ are both equal to $\frac{1}{24}\delta_A$.
Here $A$ is the graph with two genus $0$ vertices connected by $b$
edges: one has marked points indexed by $S$ and $a+1$ self-loops, and
the other has marked points indexed by $S^c$ and $c$ self-loops.

Suppose that $\gamma_1$ and $\gamma_2$ are algebraically equivalent in
$\Mla{1}{g}{n}$.  The restriction sequence
\begin{equation}\label{eq:restrict0}
  R^0(\Mbar{g}{n}\smallsetminus\Mla{i}{g}{n})\longrightarrow R^g(\Mbar{g}{n})\longrightarrow R^g(\Mla{i}{g}{n})\longrightarrow 0
\end{equation}
is exact since the first term has degree 0.  Extending $\gamma_1$ and
$\gamma_2$ to boundary strata in $\Mbar{g}{n}$,
\begin{equation}\label{eq:deltadelta}
  \gamma_1-\gamma_2\in R^0(\Mbar{g}{n}\smallsetminus\Mla{1}{g}{n}).
\end{equation}
However $R^0(\Mbar{g}{n}\smallsetminus\Mla{1}{g}{n})$ is generated by one
element $\delta_B$, where $B$ is the graph with a unique vertex of
genus $0$, $g$ self loops and $n$ half edges.

Set
\begin{align}
  {K}_1&:=2a+b+|S|   &  {L}_1&:=2c+b+|S^c|-3\\
  {K}_2&:=2c+b+|S^c|-2 &{L}_2&:= 2a+b+|S|-1.
\end{align}
For $i=1,2$ define the stratum
\begin{equation}
 \gamma_i:=\im\left(\iota_{\Gamma_i}\colon \Mbar{1}{{K}_i}\times\Mbar{0}{{L}_i+3} \to \Mbar{g}{n}\right),
\end{equation}
and
\begin{equation}
 \delta_B:=\im\left(\iota_{B}\colon \Mbar{0}{2g+n}\to\Mbar{g}{n}\right).
\end{equation}
Note that ${L}_1$ and ${L}_2$ cannot both be zero for $(g,n)\neq
(2,0)$.  If they were, then necessarily $a=0$, $b=1$, and $c=1$
(recall that $a$ is non-negative while $b$ and $c$ are positive),
which implies $(g,n)=(2,0)$.

If $L_1=0$, then $L_2\neq 0$ and we have the following equations which
follow from Lemma~\ref{le:kappakappa} and the fact that $\kappa_a$
restricted to a boundary divisor if the sum of the pull-backs of
$\kappa_a$ on each factor of the gluing map.
\begin{align}
  \kappa_{{K}_1}\gamma_1&= \frac{1}{2^{g-1}(g-2)!}\frac{1}{24} &   \kappa_1\kappa_{{L}_2}\gamma_1&= 0\\
  \kappa_{{K}_1}\gamma_2&= 0 &\kappa_1\kappa_{{L}_2}\gamma_2&=\frac{1}{2^{g-1}(g-1)!}\frac{1}{24} \\
  \kappa_{{K}_1}\delta_B&=\frac{1}{2^{g}g!}  &\kappa_1\kappa_{{L}_2}\delta_B&=\frac{1}{2^{g}g!}\left(\binom{2g+n-1}{2}-1\right)
\end{align}
These are incompatible with equation~\eqref{eq:deltadelta}.  If
$L_2=0$ and $L_1\neq 0$ a similar argument holds.

Now consider the final case where both $L_1$ and $L_2$ are nonzero.
The equations
\begin{align}
 \kappa_{{K}_1}\kappa_{{L}_1}\gamma_1&=\frac{1}{24}
  &\kappa_{{K}_1}\kappa_{{L}_1}\gamma_2&= 0
  &\kappa_{{K}_1}\kappa_{{L}_1}\delta_B&= \binom{2g+n-1}{{K}_1+1}-1 \\
\kappa_{2g-3+n}\gamma_1&= 0
  &\kappa_{2g-3+n}\gamma_2&= 0
  &\kappa_{2g-3+n}\delta_B&= \frac{1}{2^g g!}
\end{align}
show the independence of the strata $\gamma_1 - \gamma_2$ and
$\delta_B$.  Thus $\gamma_1$ cannot be algebraically equivalent to
$\gamma_2$ in $\Mla{1}{g}{n}$.
\end{proof}

\begin{corollary}\label{cor:dimension}
The dimension of the kernel of the map
\begin{equation}\label{eq:phi}
  \phi \colon R^g(\Mla{1}{g}{n})\to\hom\left(R^{2g-4+n}(\Mla{1}{g}{n}),R^{3g-4+n}(\Mla{1}{g}{n})\right).
\end{equation}
goes to infinity as $g$ or $n$ go to infinity.
\end{corollary}

\begin{proof}
  We exhibit a set of roughly $g+n/2$ linearly independent classes in $\ker \phi$. Let $\overline{n}$ denote the set
  $\{1,\ldots,n\}$. Set
  \begin{equation}
    \gamma_i = \delta(0,1,g-1,\iibar)-\delta(g-2,1,1,\nbar\smallsetminus\iibar)
  \end{equation}
  for $i=1,\ldots,\lfloor n/2\rfloor$, and
  \begin{equation}
    \eta_j = \delta(j,1,g-j-1,\nbar)-\delta(g-j-2,1,j+1,\emptyset)
  \end{equation}
  for $j=1,\ldots,g-2$.  Since
  $\lambda_1\gamma_i=\lambda_1\eta_j=0$, the classes lie in the
  kernel of $\phi$. The equations
  \begin{align}
    \psi_1^{2i+1}\kappa_{2g+n-4-2i}\gamma_j &= A_{i}\delta_{ij} \\    
    \psi_1^{2i+1}\kappa_{2g+n-4-2i}\eta_j &= 0 \\
    \psi_1^{n+2i+1}\kappa_{2g-4-2i} \gamma_j &= 0 \\
    \psi_1^{n+2i+1}\kappa_{2g-4-2i} \eta_j &= B_{i}\delta_{ij}
  \end{align}
  (here $\delta_{ij}$ is Kronecker's delta - not a boundary stratum, and $A_{i},B_{i}$ are nonzero real numbers) show that these classes are independent.  Modulo
  $R^0(\Mbar{g}{n}\smallsetminus\Mla{1}{g}{n})$, which is
  one-dimensional, roughly $g+\left\lfloor n/2\right\rfloor-1$ of them
  must remain independent.
\end{proof}

\begin{figure}[tb]
\begin{pspicture}(-1,-1)(6,1)
\rput(-1.4,0){$\gamma_1$:}
\rput(0,0){
  \genusvertex{(0,0)}{v1}{1}
  \cnode*(1,0){\nr}{v2}
  \pnode(1,0.7){x1}
  \pnode(1.7,0){x2}
  \pnode(1,-.7){x3}
  \markedpoint{(-.7,0.3)}{p1}{1}
  \markedpoint{(-.7,-.3)}{p2}{2}
  \markedpoint{(0.4,0.2)}{p3}{3}
  \markedpoint{(0.6,0.6)}{p4}{4}
  \ncline{v1}{v2}
  \ncline{v1}{p1}
  \ncline{v1}{p2}
  \ncline{v2}{p3}
  \ncline{v2}{p4}
  \nccurve[angleA=135,angleB=180]{v2}{x1}  
  \nccurve[angleA=45,angleB=0]{v2}{x1}
  \nccurve[angleA=45,angleB=90]{v2}{x2}
  \nccurve[angleA=-45,angleB=-90]{v2}{x2}  
  \nccurve[angleA=-45,angleB=0]{v2}{x3}
  \nccurve[angleA=-135,angleB=180]{v2}{x3}  
}
\rput(2,0){$-$}
\rput(3,0){
  \cnode*(0,0){\nr}{v1}
  \genusvertex{(1,0)}{v2}{1}
  \pnode(0,-.7){x1}
  \pnode(1.7,0){x2}
  \pnode(1,-.7){x3}
  \markedpoint{(-.7,0)}{p1}{1}
  \markedpoint{(0,0.7)}{p2}{2}
  \markedpoint{(0.5,0.6)}{p3}{3}
  \markedpoint{(1,0.7)}{p4}{4}
  \ncline{v1}{v2}
  \ncline{v1}{p1}
  \ncline{v1}{p2}
  \ncline{v2}{p3}
  \ncline{v2}{p4}
  \nccurve[angleA=-45,angleB=0]{v1}{x1}
  \nccurve[angleA=-135,angleB=180]{v1}{x1}  
  \nccurve[angleA=45,angleB=90]{v2}{x2}
  \nccurve[angleA=-45,angleB=-90]{v2}{x2}  
  \nccurve[angleA=-45,angleB=0]{v2}{x3}
  \nccurve[angleA=-135,angleB=180]{v2}{x3}  
}
\end{pspicture}
\begin{pspicture}(-1,-1)(6,1)
\rput(-1.4,0){$\gamma_2$:}
\rput(0,0){
  \genusvertex{(0,0)}{v1}{1}
  \cnode*(1,0){\nr}{v2}
  \pnode(1,0.7){x1}
  \pnode(1.7,0){x2}
  \pnode(1,-.7){x3}
  \markedpoint{(-.7,0.6)}{p1}{1}
  \markedpoint{(-.7,0.2)}{p2}{2}
  \markedpoint{(-.7,-.2)}{p3}{3}
  \markedpoint{(-.7,-.6)}{p4}{4}
  \ncline{v1}{v2}
  \ncline{v1}{p1}
  \ncline{v1}{p2}
  \ncline{v1}{p3}
  \ncline{v1}{p4}
  \nccurve[angleA=135,angleB=180]{v2}{x1}
  \nccurve[angleA=45,angleB=0]{v2}{x1}
  \nccurve[angleA=45,angleB=90]{v2}{x2}
  \nccurve[angleA=-45,angleB=-90]{v2}{x2}  
  \nccurve[angleA=-45,angleB=0]{v2}{x3}
  \nccurve[angleA=-135,angleB=180]{v2}{x3}  
}
\rput(2,0){$-$}
\rput(3,0){
  \cnode*(0,0){\nr}{v1}
  \genusvertex{(1,0)}{v2}{1}
  \pnode(1,0.7){x1}
  \pnode(0,-.7){x2}
  \pnode(1,-.7){x3}
  \markedpoint{(0,0.7)}{p1}{1}
  \markedpoint{(-.3,0.6)}{p2}{2}
  \markedpoint{(-.6,0.3)}{p3}{3}
  \markedpoint{(-.7,0)}{p4}{4}
  \ncline{v1}{v2}
  \ncline{v1}{p1}
  \ncline{v1}{p2}
  \ncline{v1}{p3}
  \ncline{v1}{p4}
  \nccurve[angleA=135,angleB=180]{v2}{x1}  
  \nccurve[angleA=45,angleB=0]{v2}{x1}
  \nccurve[angleA=-45,angleB=0]{v1}{x2}
  \nccurve[angleA=-135,angleB=180]{v1}{x2}  
  \nccurve[angleA=-45,angleB=0]{v2}{x3}
  \nccurve[angleA=-135,angleB=180]{v2}{x3}  
}
\end{pspicture}

\begin{pspicture}(-1,-1)(6,1)
\rput(-1.4,0){$\eta_1$:}
\rput(0,0){
  \genusvertex{(0,0)}{v1}{1}
  \cnode*(1,0){\nr}{v2}
  \pnode(1,0.7){x1}
  \pnode(0,-.7){x2}
  \pnode(1,-.7){x3}
  \markedpoint{(0,0.7)}{p1}{1}
  \markedpoint{(-.3,0.6)}{p2}{2}
  \markedpoint{(-.6,0.3)}{p3}{3}
  \markedpoint{(-.7,0)}{p4}{4}
  \ncline{v1}{v2}
  \ncline{v1}{p1}
  \ncline{v1}{p2}
  \ncline{v1}{p3}
  \ncline{v1}{p4}
  \nccurve[angleA=135,angleB=180]{v2}{x1}
  \nccurve[angleA=45,angleB=0]{v2}{x1}
  \nccurve[angleA=-45,angleB=0]{v1}{x2}
  \nccurve[angleA=-135,angleB=180]{v1}{x2}  
  \nccurve[angleA=-45,angleB=0]{v2}{x3}
  \nccurve[angleA=-135,angleB=180]{v2}{x3}  
}
\rput(2,0){$-$}
\rput(3,0){
  \cnode*(0,0){\nr}{v1}
  \genusvertex{(1,0)}{v2}{1}
  \pnode(-.7,0){x1}
  \pnode(0,-.7){x2}
  \pnode(1.7,0){x3}
  \markedpoint{(-.6,0.7)}{p1}{1}
  \markedpoint{(-.2,0.7)}{p2}{2}
  \markedpoint{(0.2,0.7)}{p3}{3}
  \markedpoint{(0.6,0.7)}{p4}{4}
 \ncline{v1}{v2}
  \ncline{v1}{p1}
  \ncline{v1}{p2}
  \ncline{v1}{p3}
  \ncline{v1}{p4}
  \nccurve[angleA=135,angleB=90]{v1}{x1}
  \nccurve[angleA=-135,angleB=-90]{v1}{x1}  
  \nccurve[angleA=-135,angleB=180]{v1}{x2}  
  \nccurve[angleA=-45,angleB=0]{v1}{x2}
  \nccurve[angleA=45,angleB=90]{v2}{x3}
  \nccurve[angleA=-45,angleB=-90]{v2}{x3}  
}
\end{pspicture}
\begin{pspicture}(-1,-1)(6,1)
\rput(-1.4,0){$\eta_2$:}
\rput(0,0){
  \genusvertex{(0,0)}{v1}{1}
  \cnode*(1,0){\nr}{v2}
  \pnode(-.7,0){x2}
  \pnode(0,-.7){x3}
  \pnode(1.7,0){x4}
  \markedpoint{(-.6,0.7)}{p1}{1}
  \markedpoint{(-.2,0.7)}{p2}{2}
  \markedpoint{(0.2,0.7)}{p3}{3}
  \markedpoint{(0.6,0.7)}{p4}{4}
  \ncline{v1}{v2}
  \ncline{v1}{p1}
  \ncline{v1}{p2}
  \ncline{v1}{p3}
  \ncline{v1}{p4}
  \nccurve[angleA=135,angleB=90]{v1}{x2}  
  \nccurve[angleA=-135,angleB=-90]{v1}{x2}
  \nccurve[angleA=-135,angleB=180]{v1}{x3}
  \nccurve[angleA=-45,angleB=0]{v1}{x3}  
  \nccurve[angleA=45,angleB=90]{v2}{x4}
  \nccurve[angleA=-45,angleB=-90]{v2}{x4}  
}
\rput(2,0){$-$}
\rput(3,0){
   \cnode*(0,0){\nr}{v1}
   \genusvertex{(1,0)}{v2}{1}
  \pnode(0,0.7){x1}
  \pnode(-.7,0){x2}
  \pnode(0,-.7){x3}
  \markedpoint{(0.5,0.6)}{p1}{1}
  \markedpoint{(0.6,0.2)}{p2}{2}
  \markedpoint{(0.6,-.2)}{p3}{3}
  \markedpoint{(0.5,-.6)}{p4}{4}
  \ncline{v1}{v2}
  \ncline{v1}{p1}
  \ncline{v1}{p2}
  \ncline{v1}{p3}
  \ncline{v1}{p4}
  \nccurve[angleA=45,angleB=0]{v1}{x1}
  \nccurve[angleA=135,angleB=180]{v1}{x1}
  \nccurve[angleA=135,angleB=90]{v1}{x2}
  \nccurve[angleA=-135,angleB=-90]{v1}{x2}
  \nccurve[angleA=-135,angleB=180]{v1}{x3}
  \nccurve[angleA=-45,angleB=00]{v1}{x3}  
}
\end{pspicture}
\caption{The classes $\gamma_1$, $\gamma_2$, $\eta_1$, and $\eta_2$ for $(g,n)=(4,4)$}
\end{figure}
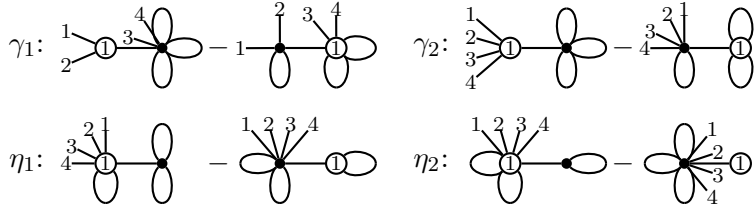

\section{Failure of Poincar\'e duality for $\Mla{i}{g}{n}$}

We first outline the strategy of proof for Proposition \ref{prop:lai}.
Fix $g$ and $i\leq g-1$. A generalization of the construction in
\S\ref{failure} produces a set $S_m\subseteq R^g(\Mbar{g}{(i+1) +m})$
in the annihilator of $\lambda_i$. The first problem in showing that
at least one such class is nonzero in $R^g(\Mla{i}{g}{(i+1) +m})$ is
that the kernel of the restriction sequence (\ref{eq:restrict}) is not
known to be tautological. This is a difficult question that we cannot
tackle at present. We therefore assume such kernel to be
tautological. Even so, the dimension of
$R^{i-1}(\Mbar{g}{n}\smallsetminus\Mla{i}{g}{n})$ grows quickly as $g$
or $n$ increase. For a fixed $g$ we bound the order of growth by
$i^n$. By proving that $S_m$ spans a linear subspace of dimension
$(i+1)^m$, we conclude that eventually some classes in $S_m$ will be
nonzero in $R^g(\Mla{i}{g}{(i+1) +m})$.

\begin{proof}[Proof of Proposition~\ref{prop:lai}] 
Let $\sigma \in R^g(\Mbar{g}{i+1})$ be the difference of boundary
classes illustrated in Figure \ref{fig:sigma}.  Intersecting either of
the two strata with $\lambda_i$ results in $1/24^i$ times the class
of the graph where all genus one vertices are replaced with loops, and
thus $\sigma\lambda_i=0$. We set $S_0=\{\sigma\}$.

\begin{figure}[tb]
\centerline{
\begin{pspicture}(0,-1)(15,1)
\rput(0,0){
  \pnode(-.7,0){x0}
  \cnode*(0,0){2pt}{v0}
  \genusvertex{(1,0)}{v1}{1}
  \genusvertex{(2,0)}{v2}{1}
  \rput(3,0){\rnode{v3}{$\cdots$}}
  \genusvertex{(4,0)}{v4}{1}
  \genusvertex{(5,0)}{v5}{1}
  \genusvertex{(6,0)}{v6}{1}
  \markedpoint{(0,0.7)}{mp0}{1}
  \markedpoint{(1,0.7)}{mp1}{2}
  \markedpoint{(2,0.7)}{mp2}{3}
  \markedpoint{(4,0.7)}{mp4}{i-1}
  \markedpoint{(5,0.7)}{mp5}{i}
  \markedpoint{(6,0.7)}{mp6}{i+1}
  \ncline{v0}{v1}
  \ncline{v1}{v2}
  \ncline{v2}{v3}
  \ncline{v3}{v4}
  \ncline{v4}{v5}
  \ncline{v5}{v6}
  \ncline{v0}{mp0}
  \ncline{v1}{mp1}
  \ncline{v2}{mp2}
  \ncline{v4}{mp4}
  \ncline{v5}{mp5}
  \ncline{v6}{mp6}
  \nccurve[angleA=135,angleB=90]{v0}{x0}
  \nccurve[angleA=-135,angleB=-90]{v0}{x0}
}
\rput(7,0){---}
\rput(8,0){
  \genusvertex{(0,0)}{v0}{1}
  \genusvertex{(1,0)}{v1}{1}
  \genusvertex{(2,0)}{v2}{1}
  \rput(3,0){\rnode{v3}{$\cdots$}}
  \genusvertex{(4,0)}{v4}{1}
  \genusvertex{(5,0)}{v5}{1}
  \cnode*(6,0){2pt}{v6}
  \pnode(6.7,0){x6}
  \markedpoint{(0,0.7)}{mp0}{1}
  \markedpoint{(1,0.7)}{mp1}{2}
  \markedpoint{(2,0.7)}{mp2}{3}
  \markedpoint{(4,0.7)}{mp4}{i-1}
  \markedpoint{(5,0.7)}{mp5}{i}
  \markedpoint{(6,0.7)}{mp6}{i+1}
  \ncline{v0}{v1}
  \ncline{v1}{v2}
  \ncline{v2}{v3}
  \ncline{v3}{v4}
  \ncline{v4}{v5}
  \ncline{v5}{v6}
  \ncline{v0}{mp0}
  \ncline{v1}{mp1}
  \ncline{v2}{mp2}
  \ncline{v4}{mp4}
  \ncline{v5}{mp5}
  \ncline{v6}{mp6}
  \nccurve[angleA=45,angleB=90]{v6}{x6}
  \nccurve[angleA=-45,angleB=-90]{v6}{x6}
}
\end{pspicture}
}
\caption{The class $\sigma\in R^g(\Mla{i}{g}{i+1})$}\label{fig:sigma}
\end{figure}


For $m>0$, and $\mathbf{a}=(a_1,\ldots, a_m)$ an $m$-tuple of numbers
between $1$ and $i+1$, let $\sigma_\mathbf{a} \in R^g(\Mbar{g}{(i+1)
  +m})$ be the class obtained by decorating both graphs of $\sigma$
with the $j$-th mark on the $a_j$-th vertex for $j=1,\ldots,m$.  The
set $S_m$ of all possible such classes in $\Mla{i}{g}{(i+1)+m}$ has
cardinality $(i+1)^m$.

We construct inductively a set $T^m$ of classes in complementary
codimension which is dual to $S_m$.
Our base case is $m=0$, where the vector $\tau$ can be chosen to be a
scalar multiple of an appropriate product of $\psi$
classes. Consider the universal family
\begin{equation}
  \pi\colon\Mbar{g}{(i+1)+m+1} \rightarrow \Mbar{g}{(i+1)+m},
\end{equation}
and note that
\begin{equation}\label{eq:sum}
  \pi^\ast \sigma_\mathbf{a}= \sum_{k=1}^{i+1} \sigma_{(\mathbf{a},k)}.
\end{equation}
where $(\mathbf{a},k)$ is the sequence with $k$ appended to the end of
$\mathbf{a}$.  If $D_{j,m+1}$ denotes the divisor image of the $j$-th
section in $\Mbar{g}{(i+1)+m+1}$, then for any $m$-tuple $\mathbf{a}$,
\begin{equation}
\label{kill}
D_{j,m+1} \sigma_{(\mathbf{a},k)}=0
\end{equation}
if $j \not = k$.  By the projection formula, and equations
(\ref{eq:sum}) and (\ref{kill}),
\begin{equation}
D_{j,m+1}\pi^\ast(\tau^\mathbf{b}) \sigma_{(\mathbf{a},k)}
\end{equation}
equals the class of a point if $j=k$ and $\mathbf{a}=\mathbf{b}$, and
vanishes otherwise. Therefore the set:
\begin{equation}
T^{m+1}:=\{\pi^\ast\tau^{b_1,\ldots,b_m}D_{b_{m+1},m+1}\}_\mathbf{b} 
\end{equation}
gives a dual basis to $S_{m+1}$.

The growth of $\dim R^{i-1}(\Mbar{g}{n}\smallsetminus\Mla{i}{g}{n})$
with respect to $n$ is at most $O(i^n)$.  To see this, note that the
decorated dual graph of any class in
$R^{i-1}(\Mbar{g}{n}\smallsetminus\Mla{i}{g}{n})$ has at most $i$
vertices, and $e$ edges, where $g-i+1\leq e\leq g$.  The total number
of possibly unstable graphs without marked points satisfying these
conditions is independent of $n$.

Graphs with exactly $i$ vertices have $g$ edges, and hence classes
supported on these strata are pure boundary. For a given graph there
are $i^n$ possible ways of distributing the marked points on the
vertices.

For graphs with strictly less than $i$ vertices, there are at most
$(i-1)^n$ ways to distribute the marks on the vertices.  Each vertex
can be decorated with a monomial in $\psi$ and $\kappa$
classes~\cite{MR1960923}*{Proposition 11}, of degree $\leq
i-1$~\cite{MR1960923}*{Proposition 11}.  The number of ways to choose
$\kappa$ classes to decorate the vertices is independent of $n$.  The
number of monomials in $\psi$ classes of degree $\leq i-1$ grows
polynomially in $n$, yielding an order of $O((i-1)^n n^k)<O(i^n)$ for
the number of classes supported on strata with less than $i$ vertices.

Altogether, we have obtained that the dimension of
$R^{i-1}(\Mbar{g}{n}\smallsetminus\Mla{i}{g}{n})$ grows at most as
$O(i^n)$. Thus $\dim R^g(\Mla{i}{g}{n}) > \dim
R^{i-1}(\Mbar{g}{n}\smallsetminus\Mla{i}{g}{n})$ for large $n$, which
implies that some classes in the annihilator of $\lambda_i$ are
nonzero in $R^g(\Mla{i}{g}{n})$.
\end{proof}

All nontrivial classes in the annihilators of $\lambda_i$ that we
construct are in codimension $g$ or higher.  This leaves us with a
natural question:

\begin{question}
  Is the map
\begin{equation}
  \phi \colon R^j(\Mla{i}{g}{n})\to\hom\left(R^{3g-3+n-i-j}(\Mla{i}{g}{n}),R^{3g-3+n-i}(\Mla{i}{g}{n})\right)
\end{equation}
  injective for $j\leq g-1$?
\end{question}

\end{document}